\documentclass[12pt]{article}
\usepackage{amssymb,amsmath,graphicx,hyperref}
\usepackage[utf8]{inputenc}
\usepackage[english]{babel}

\title{Counting Truchet Tile Balls}
\author{Thomas Fernique}
\date{thomas.fernique@ens-lyon.org}

\begin{document}
\maketitle

\begin{abstract}
A formula is established that counts the number of different balls that can be made by decorating the pentagons and hexagons of a classic football with Truchet-like patterns.
\end{abstract}

\begin{figure}[hbtp]
\centering
\includegraphics[height=28mm]{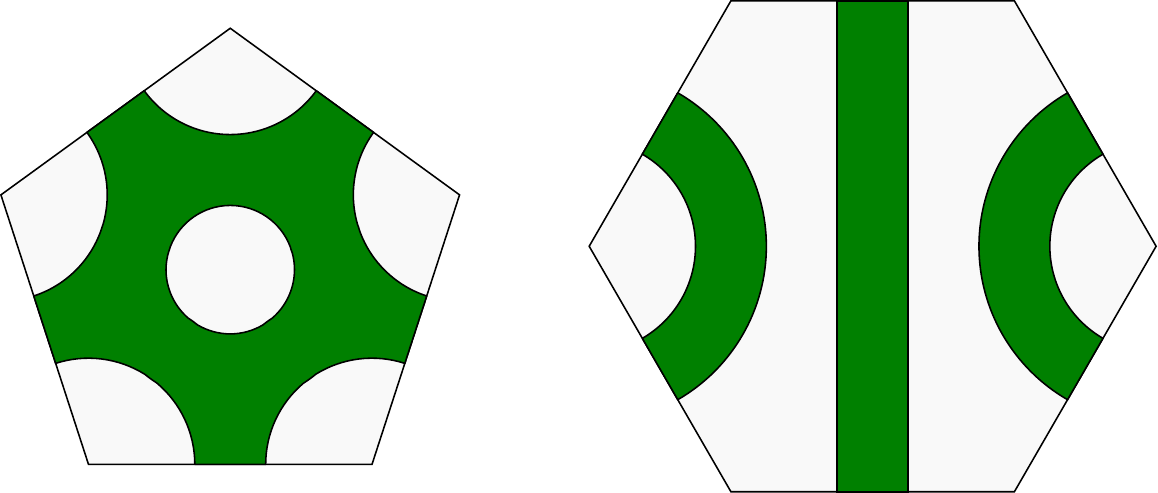}
\caption{The Truchet-like panels designed by Jon-Paul Wheatley \cite{Whe24}.}
\label{fig:truchet_base}
\end{figure}

Among the numerous ball designs created by Jon-Paul Wheatley, one is based on two panels looking like the ones depicted in Figure~\ref{fig:truchet_base}.
These panels are stitched to form a truncated icosahedron ($12$ pentagons and $20$ hexagons), as in a classic football (the Telstar ball).

\begin{figure}[hbtp]
\centering
\includegraphics[width=\textwidth]{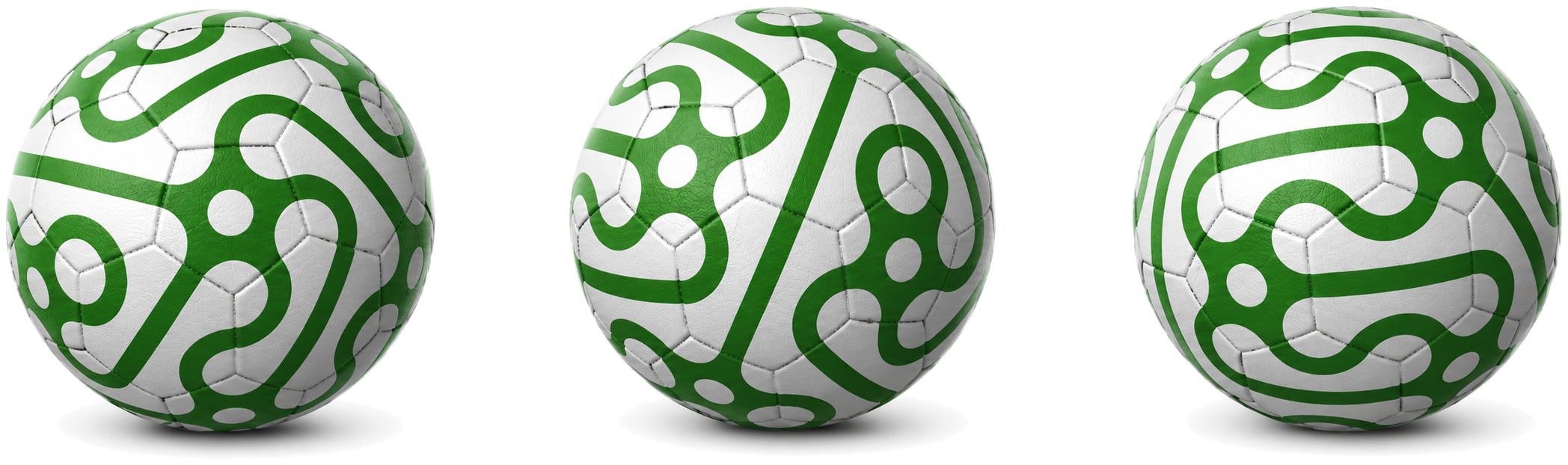}
\caption{Some footballs made with the panels of Fig.~\ref{fig:truchet_base} (rendering).}
\label{fig:balls_base}
\end{figure}

The point is that each panel can be randomly oriented, independently of its neighbors.
Indeed, two neighbor panels always match because each side sees the same pattern, namely a centered thick colored line.
This allows one to make a lot of different designs, see Figure~\ref{fig:balls_base}.
It is reminiscent of the square tile introduced by Sébastien Truchet in 1704, even more the variant introduced by Smith in 1987 \cite{SB87}.
This design also appeared in several board games, e.g., {\em Black Path Game} (1960), {\em Trax} (1980) or {\em Tantrix} (1991).

\pagebreak[1]

We here address the question of the number of different balls that can be made.
A crude estimate is $3^{20}$, that is, roughly 3.48 billion, because each of the 20 hexagons can appear in three different orientations (the orientation of pentagons does not play a role because of the symmetry of its design).
However, many of these balls are actually identical up to a rigid motion.
Can we find the exact number of different balls?

The key tool is the Burnside lemma \cite{Bur87}.
Namely, given a finite group $G$ that acts on a set $X$, the number of orbits $|X/G|$ satisfies:
\begin{equation}
\label{eq:burnside}
|X/G|=\frac{1}{|G|}\sum_{g\in G}|X^g|,
\end{equation}
where $|X^g|$ denotes the number of elements in $X$ that are fixed by $g$.
In our case, $G$ is the set of rigid motions that fix the truncated icosahedron, $X$ is the set of balls that can be formed with our Truchet panels without considering rigid motions, and $|X/G|$ is the size of this set up to rigid motion, that is, the number we are interested in.
Hence, all we have to do is to find the size $|G|$ of $G$ and, for each $g$ in $G$, the number $|X^g|$ of Truchet tile balls that are fixed by the rigid motion $g$.

The $60$ rigid motions that fix the truncated icosahedron are well known:
\begin{enumerate}
\item
The rotations of angle $\pm \tfrac{2\pi}{5}$ ($\pm 72^\circ$) and $\pm \tfrac{4\pi}{5}$ ($\pm 144^\circ$) around the line that joins the centers of two opposite pentagons.
There are 6 such lines, thus the four above angles yield 24 rotations.
\item
The rotations of angle $\pm \tfrac{2\pi}{3}$ ($\pm 120^\circ$) around the line that joins the centers of two opposite hexagons.
There are 10 such lines, thus 20 such rotations.
\item
The rotation of angle $\pi$ ($180^\circ$) around the line that joins the centers of two opposite edges shared by two hexagons.
There are 15 such lines, thus 15 such rotations.
\item
The identity, that is, the rigid motion that consists of doing nothing.
\end{enumerate}

\pagebreak[1]

\begin{figure}[hbtp]
\centering
\includegraphics[height=28mm]{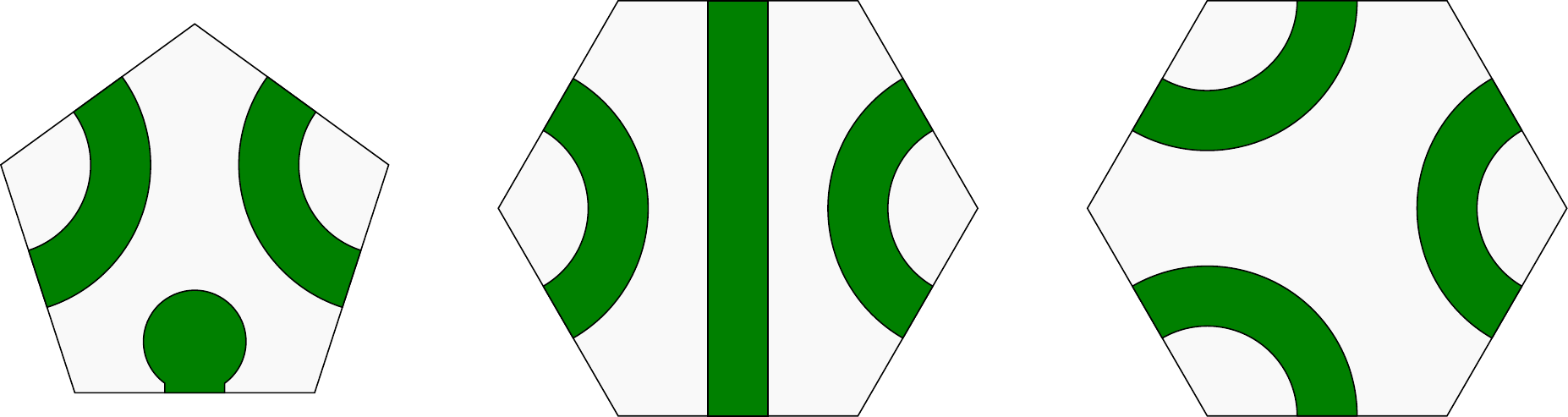}
\caption{Another set of Truchet-like panels.}
\label{fig:truchet_rab}
\end{figure}

What about $|X^g|$?
For the sake of generality, let us denote by $p$ and $q$ the numbers of different patterns that can appear on, respectively, a pentagonal and a hexagonal panel (rotated patterns are considered to be different).
Fig.~\ref{fig:truchet_base} corresponds to $p=1$ and $q=3$, but other patterns could be considered (Fig.~\ref{fig:truchet_rab}).
We make cases, depending on $g\in G$:
\begin{enumerate}
\item
If $g$ is a rotation of angle $\pm\frac{2\pi}{5}$ or $\pm\frac{4\pi}{5}$, then every face of the truncated icosahedron -- but the two pentagons crossed by the axis of the rotation -- has an orbit of size 5, that is, applying $g$ over and over on this face yields exactly 5 different faces.
The two pentagons crossed by the axis have an orbit of size 1.
There is thus $(12-2)/5+2=4$ orbits of pentagons and $20/5=4$ orbits of hexagons.
We can freely choose the pattern drawn on exactly one panel in each orbit: the patterns on the other panels are then obtained by applying the rotation.
This yields $|X^g|=p^4q^4$.
\item
If $g$ is a rotation of angle $\pm\frac{2\pi}{3}$, then every face of the truncated icosahedron -- but the two hexagons crossed by the axis of the rotation -- has an orbit of size 3.
The two hexagons crossed by the axis have an orbit of size 1.
There is thus $12/3=4$ orbits of pentagons and $(20-2)/3+2=8$ orbits of hexagons.
This yields $|X^g|=p^4q^8$.
\item
If $g$ is a rotation of angle $\pi$, then every face has an orbit of size 2.
There is thus $12/2=6$ orbits of pentagons and $20/2=10$ orbits of hexagons.
This yields $|X^g|=p^6q^{10}$.
\item
If $g$ is the identity, then every face has an orbit of size 1.
There is thus $12$ orbits of pentagons and $20$ orbits of hexagons.
This yields $|X^g|=p^{12}q^{20}$.
\end{enumerate}
We can now plug in Equation~\eqref{eq:burnside} the number $|X^g|$ for every type of rigid motion $g\in G$ and use our count of every type of rigid motion in $G$ to get the total number of different balls up to a rigid motion:
\begin{equation}
\label{eq:main}
|X/G|=\tfrac{1}{60}(24p^4q^4+20p^4q^8+15p^6q^{10}+p^{12}q^{20}).
\end{equation}

\pagebreak[1]

For $p=1$ and $q=3$, that is, the patterns designed by Jon-Paul Wheatley and depicted in Fig.~\ref{fig:truchet_base}, this yields $58\,130\,055$ different balls.
This is more or less the number of footballs manufactured every year in the whole world.
This is also almost exactly $60$ times less than the crude estimate $3^{20}$, what is not that surprising because random patterns most of the time break all the symmetries of the truncated icosahedron, hence the number of balls up to a rigid motion is roughly divided by $|G|=60$.
Putting all these balls on the ground as densely as possible (that is, on a triangular lattice) will cover around 243 hectares -- roughly the area of Hyde Park in London.
Stacking them as densely as possible (as oranges on a market stall) will fill a cube with 76-meter sides or 175 Olympic-sized swimming pools.

If we use, instead, the Truchet-like panels depicted in Fig.~\ref{fig:truchet_rab}, then $p=5$ (the five orientations of the pentagon are different) and $q=3+2=5$ (3 for the first hexagon, 2 for the second one).
Equation~\eqref{eq:main} yields that there are $388\,051\,072\,794\,677\,890\,625$ different balls!
Stacked as densely as possible, they would cover the entire Earth, including the oceans (assuming the balls float), to a depth of about $5\,700$ meters — higher than Mount Ararat ($5\,137$m).
A veritable Flood of balls!

\begin{figure}[hbtp]
\centering
\includegraphics[width=\textwidth]{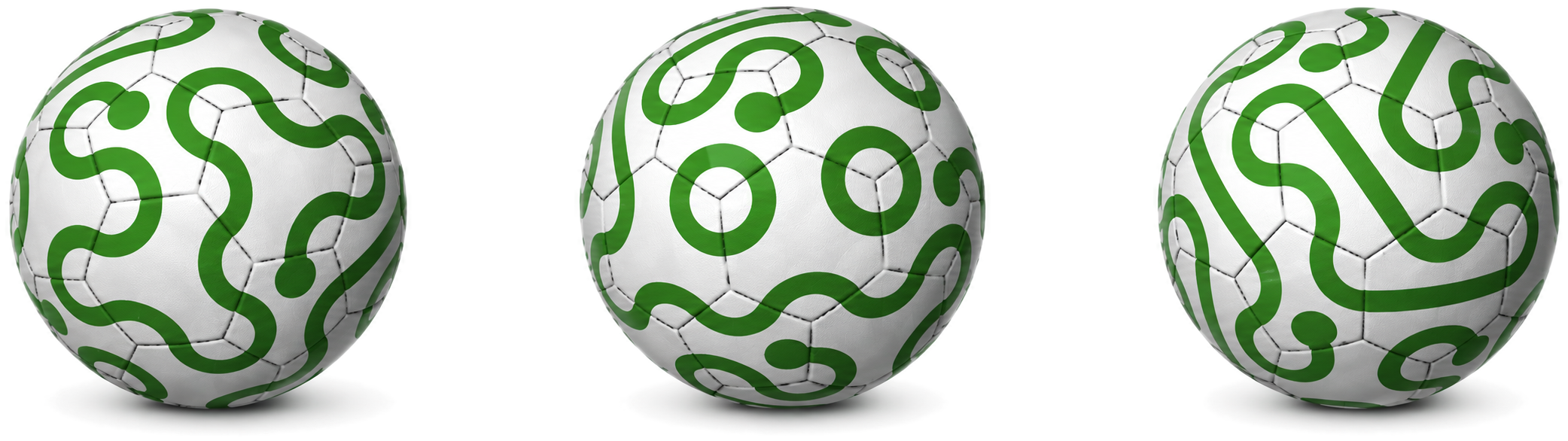}
\caption{Some footballs made with the panels of Fig.~\ref{fig:truchet_rab} (rendering).}
\label{fig:balls_rab}
\end{figure}

\bibliographystyle{alpha}
\bibliography{truchet}

\end{document}